\documentclass[a4paper]{amsart}
\usepackage{amsmath}
\def\me{\mathbf{g}}
\def\contr{\neg}

\def\ka{K\"{a}hler~}
\def\ak{affine-\ka}
\def\ma{Monge-Amp\`{e}re\,\,}
\def\contr{\rightharpoonup}

\def\rnum{\mathbb{R}}
\def\cnum{\mathbb{C}}
\def\znum{\mathbb{Z}}
\def\zint{\mathbb{Z}}
\def\ka{K\"{a}hler\;}

\def\sud{{\mathcal{S}^{1\to2}}}
\def\sdu{{\mathcal{S}^{2\to1}}}
\def\R{\mathbb{R}}

\def\ei{{E^{0}_i}}

\def\eik{{E^{k}_{i}}}

\def\ii{{E^{\bar{0}}_i}}
\def\iik{{E^{\bar{k}}_{i}}}

\def\eial{{E^\alpha_i}}
\def\ejbe{{E^\beta_j}}
\def\eibe{{E^\beta_i}}
\def\eibaral{{E^{\bar{\alpha}}_i}}
\def\eibarbe{{E^{\bar{\beta}}_i}}
\def\ejga{{E^\gamma_j}}
\def\eiga{{E^\gamma_i}}
\def\eibarga{{E^{\bar{\gamma}}_i}}
\def\ejbarga{{E^{\bar{\gamma}}_j}}
\def\ejbaral{{E^{\bar{\alpha}}_j}}
\def\ejbarbe{{E^{\bar{\beta}}_j}}
\def\lop{\mathbf{L}}
\def\x12{\mathbf{X_1}\times_B \mathbf{X_2}}
\def\mex12{\me_{_{\x12}}}

\def\me{\mathbf{g}}

\def\contraction{\lrcorner}
\newtheorem{teo}{Theorem}[section]
\newtheorem{con}{Conjecture}[section]
\newtheorem{cor}[teo]{Corollary}
\newtheorem{exe}[teo]{Example}
\newtheorem{lem}[teo]{Lemma}
\newtheorem{pro}[teo]{Proposition}
\newtheorem{dfn}[teo]{Definition}

\newtheorem{rmk}[teo]{Remark}

\title{Mirror symmetry and self-dual manifolds}
\author{Michele Grassi}
\date{February 2, 2002}
\begin{document}
\begin{abstract}
We introduce self-dual manifolds and show that they can be used to encode mirror symmetry for \ak manifolds and for elliptic curves. Their geometric properties, especially the link with special lagrangian fibrations and the existence of a transformation similar to the Fourier-Mukai functor, suggest that this approach may be able to explain mirror symmetry also in other situations.
\end{abstract}
\maketitle
\section{Introduction}
Starting with the paper ~\cite{SYZ} by Strominger, Yau and Zaslow, the search for a geometric counterpart to  mirror symmetry  has beed directed mainly at special lagrangian fibrations of Calabi-Yau manifolds. In recent years however this approach has come under some criticism, because it appears unlikely that such fibrations will exist in general, at least of the well behaved kind required for mirror symmetry. Much of the research has therefore usually assumed the existence of such well behaved fibrations (see for example ~\cite{Gross}), and  has studied the behaviour ''in the large complex structure limit''  (\cite{GW}, ~\cite{KS}).
In this paper we propose that the "right" object to associate to a mirror pair of Calabi-Yau manifolds should not be a pair of lagrangian fibrations, but a self-dual manifold. As will become clear later, the special lagrangian fibrations then come back into the picture as a "Gromov-Hausdorff limits" of foliations on the self-dual manifold, and as such may be very singular and badly behaved in general. Instead, the self-dual manifold is expected to be smooth, and to contain in its structure the tools to explain mirror symmetry.
We prove that this picture is correct in the case for complex dimension $n=1$ (elliptic curves), and for \ak manifolds of any dimension. Although we did not attempt to include them here for reasons of space, we believe that the cases of abelian varieties and of $K3$'s should be within the reach of the tecniques that we introduce.\\ 
The basic tool the following geometric data: a smooth Riemannian manifold $(X,\me)$ of (real) dimension $3n$, together with two smooth differential forms of degree $2$ on $X$, $\omega_1$ and $\omega_2$, which are compatible  with the metric, in the sense that at all points $p\in X$ there is an orthonormal coframe 
\[dx_1,...,dx_n,dy^1_1,...,dy^1_n,dy^2_1,...,dy^2_n\]
such that $\omega_j = \sum_i dx_i\wedge dy^j_i$. Self-dual manifolds are objects of the above kind, with two more conditions on the data.  For the precise definition, see the next section. \\
Of course, the final goal of explaining mirror symmetry in terms of self-dual manifolds, even if achievable, will require a lot of effort. In the present paper we provide some clues as to why we think our approach should work.\\
A self-dual manifold of real dimension $3n$ should be a way to interpolate between two mirror dual manifolds of complex dimension $n$, which can be recovered as natural Gromov-Hausorff limits. Again, we prove this statement only for \ak manifolds and for elliptic curves. In general, this interpolation property is expected to happen at limit points in the moduli spaces of the mirror pair. However, in the case of elliptic curves and of \ak manifolds we show that this holds at all points.\\
One of the advantages of self-dual manifolds over the traditional approach via special lagrangian fibrations is that while in general the fibrations are expected to exist only in the limit, you just expect a self-dual manifold also at finite points; the original Calabi-Yaus are then quotients (with respect to foliations), which near the boundary of the moduli spaces became Gromov-Hausdorff limits.\\ 
To build the smooth self-dual manifold associated to a mirror pair we start  from the fibre product of dual special lagrangian fibrations, which in our case don't have singular fibres. It should be possible to use this method also when the fibrations have isolated singularities. We did not elaborate on this in the present paper.\\
Another advantage of self-dual manifolds is that their structure can be significantly weakened (to a polysymplectic structure) or strengthened (to a $2$-\ka structure). While polysymplectic manifolds share with symplectic ones the absence of local moduli, $2$-\ka manifolds are in a sense similar to hyperk\"{a}hler ones (although they have dimension 3n). 
The rich algebraic structure of the cohomology of $2$-\ka manifolds is what brought us to their study in the first place, although we were not very successfull in constructing smooth compact ones (except in the homogeneous cases). This might be just a temporary limitation, or might be due to some actual obstructions. In any case, we expect that the  $\mathbf{sl}(4,\rnum)$  representation which exists on the cohomology of compact $2$-\ka manifolds will be useful for the study of the cohomology of self-dual manifolds near boundary points of their moduli spaces.\\
An aspect which we did not develop in the present paper is the connection of self-dual manifolds with other constructions unrelated to mirror symmetry. For example, contact structures and Seifert fibrations come into play when studying self-dual manifolds of dimension $3$. In a future paper we plan to investigate the relationship of $3n$ dimensional self-dual manifolds with $c = 3n$ (super) conformal field theory and the geometry of PDE's with target an $n$-dimensional manifold. On this last subject some material can already be found in the first part of ~\cite{G}.\\
Let us now come to a brief description of the contents of the various sections:\\
In section ~\ref{sec:a2kandsd} we introduce our main object of study, self-dual manifolds. To do that, we choose to introduce first the weaker notions of polysymplectic manifold and of almost $2$-\ka manifold, because they will play a role later. Almost $s$-\ka manifolds are polysymplectic manifolds together with a compatible metric. We show that in the case $s=2$ you have a natural {\em dualizing form}. When this form is closed, and the leaves of a certain foliation have all Riemannian volume one, you have {\em self-dual manifolds}. We show some natural ways of deforming almost $2$-\ka (and self-dual) manifolds, and finally we  introduce a transformation which is similar in nature to the Fourier-Mukai transform, and will play a role if one will want to use self-dual manifolds to prove mirror symmetry.\\
In section ~\ref{sec:slag} we show that fibre products of  Riemannian lagrangian fibrations of almost \ka manifolds over the same base give rise to almost $2$-\ka manifolds. We show that this applies to the significant case of special lagrangian tori fibrations of Calabi-Yau manifolds. Of course, you do not expect to obtain self-duality unless you start from a mirror pair.\\
In section ~\ref{sec:maandell} we apply the notions developed in the previous sections to show that self-dual manifolds can indeed be used to characterize mirror symmetry for \ak manifolds. \\
In section ~\ref{sec:ellcur} we do the same that we did in the previous section, this time for elliptic curves. We also formulate a conjecture which deals with what to expect in the case of $K3$ surfaces.\\
In section ~\ref{sec:polyandska} we show that polysymplectic manifolds have no local moduli, and we prove that the space of metrics compatible with a given polysymplectic structure is non-empty and contractible. We then introduce $2$-\ka manifolds, and prove a characterization of them which generalizes a classical one for \ka ones. We prove that a $2$-\ka manifold automatically has one of the two properties of self-dual manifolds (namely the dualizing form is closed), and that the natural deformations introduced in section ~\ref{sec:a2kandsd} carry over to the $2$-\ka case.\\
In section ~\ref{sec:reps} we show that there is a natural action of the Lie algebra $\mathbf{sl}(4,\rnum)$ on the forms of an almost $2$-\ka manifold, which in the $2$-\ka case carries over to an action on the harmonic forms.\\
\textbf{Notation}\\
For a form $\alpha$, we indicate with $\alpha^\perp$ the space of vectors which contract to zero with it. We indicate with $T^*M$ the cotangent bundle to the manifold $M$ and its total space. If $\Gamma$ is a lattice inside an euclidean space, its dual lattice (with respect to the metric) is indicated with $\Gamma^\vee$.
\section{Self-dual manifolds}
\label{sec:a2kandsd}
While the crucial notion is that of self-dual manifold (Definition ~\ref{dfnsdm}), we first introduce polysymplectic and almost $s$-\ka manifolds, which generalize symplectic and almost \ka manifolds respectively.
\begin{dfn}[\cite{G}, Definition 2.1 page 12] ~  \\
1) A {\em polysymplectic structure} on a vector space $V$  of dimension $n(s+1)$ is given by $s$ elements $\omega_1,...,\omega_s$ of $\Lambda^2V^*$ which, in some basis $v_1,...,v_n,w^1_1,...,w^s_n$  for $V$, have the {\em polysymplectic normal form}
\[\omega_j = \sum_{i=i}^nv_i^*\wedge (w^j_i)^*\]
Any such basis is called {\em standard polysymplectic}.\\
2) A {\em polysymplectic} manifold is given by a smooth manifold $X$ of dimension $n(s+1)$, together with $s$ smooth differential forms $\omega_1,...,\omega_s$ of rank $2$ such that:\\
a) The forms $\omega_j$ are closed,\\
b) At all points $p\in X$ the forms $(\omega_1)_p,...,(\omega_s)_p$ determine a polysymplectic structure on $T_pX$,\\
c) The distribution $\sum_j\omega_j^\perp$ is integrable.
\end{dfn}
The notion of polysymplectic manifold reduces to that of symplectic manifold for $s = 1$, and in that case condition $2c)$ is automatically true. The case relevant for mirror symmetry is $s = 2$, and in this case condition $2c)$ does not follow from the other ones.
\begin{exe}
\label{execotangent}
Let $M$ be a smooth manifold, and let $\mathbf{T}^{*}M$ indicate the cotangent bundle of $M$. If 
\[X ~=~\mathbf{T}^{*}M~ ~\times_{_{M}} ~~\cdots~ ~\times_{_{M}} ~~\mathbf{T}^{*}M~(s~times),\]
$\pi_{i}~:~X\to \mathbf{T}^{*}M$
is the projection on the $i^{th}$ factor, and $\omega$ the canonical symplectic form on 
$\mathbf{T}^{*}M$, let
$\omega_{i}~:=~ \pi_i^{*}\omega$.
We have then that $(X,\omega_{1},...,\omega_{s})$ is polysymplectic.
\end{exe}
\begin{dfn}[\cite{G}, Definition 6.1 page 30]
Let $(X,\omega_1,...,\omega_s)$ be a polysymplectic manifold, and let $\me$ be a Riemannian metric on $X$. We say that $\me$ is {\em compatible with the polysymplectic structure} if for every $p\in X$ there exists a standard polysymplectic basis of $T_pX$ which is also orthonormal with respect to $\me$. In that case, we say that $(X,\omega_1,...,\omega_s,\me)$ is an {\em almost $s$-\ka manifold}. 
\end{dfn}
Again, for $s = 1$ the previous notion reduces to the classical one of almost \ka manifold.
\begin{dfn}
Let $\mathbf{X} = (X,\omega_1,\omega_2,\me)$ be an almost $2$-\ka manifold. The dualizing form $\omega_D$ is the differential form of degree $2$ defined at the point $p\in X$ as
\[\sum_{i=1}^n(w_i^1)^*\wedge (w_i^2)^*\]
for any orthonormal standard polysymplectic basis $v_1,...,v_n,w_1^1,...,w_n^2$ on $T_pX$.
\end{dfn}
\begin{rmk}
\label{rmkomegawelldefined}
The form $\omega_D$ is well defined, as if $\tilde{v}_1,...,\tilde{v}_n,\tilde{w}_1^1,...,\tilde{w}_n^2$ is another orthonormal standard polysymplectic basis, it is easy to see that we must have 
\[(\tilde{w}^1_i)^* =\sum_ja_{ij}(w^1_j)^*,~(\tilde{w}^2_i)^* =\sum_ja_{ij}(w^2_j)^*\]
for some orthogonal matrix $(a_{ij})$, and therefore
\[\sum_{i}(\tilde{w}_i^1)^*\wedge \tilde{w}_i^2 = \sum_{i,j,k}a_{ij}a_{ik}(w_j^1)^*\wedge (w_k^2)^* = \sum_i(w_i^1)^*\wedge (w_i^2)^*\]
\end{rmk}
\begin{dfn}
\label{dfnsdm}
Let $\mathbf{X} = (X,\omega_1,\omega_2,\me)$ be an almost $2$-\ka manifold. We say that $\mathbf{X}$ is a {\em self-dual almost $2$-\ka manifold} (ore more briefly a {\em self-dual manifold}) if:\\
1) The differential form   $\omega_D$  is closed.\\
2) The leaves of the foliation $\omega_1^\perp + \omega_2^\perp$ have all (Riemannian) volume equal to one.\\
If only condition one holds,  $\mathbf{X} = (X,\omega_1,\omega_2,\me)$ is a  {\em weakly self-dual manifold}
\end{dfn}
\begin{exe}
Let $l \in\rnum^+$, and let
\[X = \rnum/l\znum\times \rnum/\znum\times \rnum/m\znum\]
Call $y^1$ (resp. $x$, $y^2$) the coordinate induced by $\rnum$ on the first (resp. second, third) factor. With this choice, and with  
\[\me = (dx)^2 + (dy^1)^2 + (dy^2)^2, ~\omega_1 = dx\wedge dy^1,~ \omega_2 = dx\wedge dy^2\]
we have that $(X,\omega_1,\omega_2,\me)$ is weakly self-dual. If $lm = 1$ it is also self-dual.
\end{exe}
Weakly self-dual manifolds have a very rich structure, and as we will see in the following, are rather easy to construct if you do not insist on them being compact. For now, let us just point out a feature which may look like a hyperk\"{a}hler property (although, as the above example shows, there are compact self-dual manifolds of dimension $3$)

\begin{rmk}
\label{rotatesd}
Let $\mathbf{X} = (X,\omega_1,\omega_2,\me)$ be a weakly self-dual manifold. Then, if $\omega_D$ is its dualizing form and the distribution $\omega_D^\perp + \omega_1^\perp$ is integrable, also  $(X,\omega_D,\omega_1,\me)$ is a weakly self-dual manifold. If the distribution $\omega_D^\perp + \omega_2^\perp$ is integrable, also  $(X,\omega_2,\omega_D,\me)$ is a weakly self-dual manifold.
\end{rmk}
In all the examples that we will build in this paper, both the integrability conditions of the previous remark hold. In those cases, there are three different structures of weakly self-dual manifold on the same underlying Riemannian manifold.
\begin{pro}
\label{covtosd}
Let $\nabla$ be the Levi-Civita connection associated to the metric. If $\nabla\omega_1 = \nabla\omega_2 = 0$, then $\nabla\omega_D = 0$ and hence the manifold is weakly self-dual.
\end{pro}

\textit{Proof}
If the two forms $\omega_1,\omega_2$ are covariant constant, then parallel transport, which is also orthogonal, will send any orthonormal standard polysymplectic basis to an orthonormal standard polysymplectic basis. From this and the same reasoning of Remark ~\ref{rmkomegawelldefined}, we conclude that $\omega_D$ is sent to itself.
\qed

There are various ways of deforming an almost $2$-\ka manifold. In the following definition we single out three of the most relevant ones for questions concerning Mirror Symmetry.
\begin{dfn}
Let $\mathbf{X} = (X,\omega_1,\omega_2,\me)$ be an almost $2$-\ka manifold, and let $t\in\R^+$. We then define:\\ 
1) 
\[\alpha_t\left(\mathbf{X}\right) = (X,t\omega_1,\omega_2,\alpha_t(\me))\]
where $\alpha_t(\me))$ is such that given any orthonormal standard polysymplectic basis $v_1,...,v_n$, $w^1_1,...,w^2_n$ of $T_pX$ (with respect to the given almost $2$-\ka structure on $X$), it assigns lenght squared $t$ to all the $v_i,w^1_i$ and lenght squared $t^{-1}$ to the $w^2_i$, for $i=1,...,n$.\\
2) 
\[\beta_t\left(\mathbf{X}\right) = (X,\omega_1,t\omega_2,\beta_t(\me))\]
is defined in the same way as $\alpha_t$, only with the indices $1$ and $2$ in the definition of the metric exchanged.\\
3) 
\[\lambda_t\left(\mathbf{X}\right) = (X,t\omega_1,t\omega_2,t\me)\]
\end{dfn}
We omit the easy proof of the following proposition
\begin{pro}
Let $\mathbf{X}$ be an almost $2$-\ka manifold. Then\\
1) $\alpha_t(\mathbf{X}),\beta_t(\mathbf{X}),\lambda_t(\mathbf{X})$ are almost $2$-\ka manifolds.\\
2) The deformations $\alpha_t$ and $\beta_t$ leave the dualizing form $\omega_D$ unchanged. In particular, if $\mathbf{X}$ is   (weakly) self-dual, then also $\alpha_t(\mathbf{X}),\beta_t(\mathbf{X})$ alre (weakly) self-dual.\\
3) If $\omega_i$ is (Levo-Civita) covariant constant with respect to $\me$, then it is so also with respect to $\alpha_t(\me),\beta_t(\me)$
\end{pro}
The following maps are one of the ingredients of the Mirror correspondence, and are similar in nature and in behaviour to the Fourier-Mukai functor of ~\cite{Muk}. 
\begin{dfn}
Let $(X,\omega_1,\omega_2,\me)$ be a self-dual manifold, and assume that there are surjections $\pi_1:X\to X_1$,$\pi_2:X\to X_2$ and $\pi_B:X\to B$ with compact fibres equal to leaves of the foliations $\omega_1^\perp$,$\omega_2^\perp$ and $\omega_1^\perp + \omega_2^\perp$ respectively. Let $\alpha\in\Omega^i(X_1)$. Define $\sud^{,j}(\alpha)\in \Omega^{i+2j-n}(X_2)$ as
\[\sud^{,j}(\alpha)_p(v_1,...,v_{i+2j-n}) = 
\int_{\pi_2^{-1}(p)}
\tilde{v}_1\contraction\cdots \tilde{v}_{i+2j-n}\contraction
(\omega_D^j\wedge\pi_1^*\alpha)\]
for $p\in X_2$, $v_k\in T_pX_2$ and $\tilde{v}_k$ the lifting of $v_k$ to a vector field along $\pi_2^{-1}(p)$ which projects to $v_k$ and is orthogonal to the fibre.
Define also $\sud = \sum_j \sud^{,j}$. Define $\sdu^{,j}$ and $\sdu$ similarly, but with the indices $1$ and $2$ interchanged.
\end{dfn}
Althoug the previous definition has many similarities with that of the Fourier-Mukai functor, note that we are not assuming that the fibres of $X\to B$ are tori, or flat with the induced metric. We are not using an almost complex structure on the fibres, and the natural one induced by $\omega_D$ and by the metric is "wrong" (in the case of elliptic fibrations of $K3$'s one would say that it is "rotated" with respect to the one in which the Poincar\'{e} bundle is defined). Also note that we did not attempt to find the right sign in the definition of $\sud$.
\section{Special lagrangian fibrations}
\label{sec:slag}
We now build examples of almost $2$-\ka manifolds starting from Riemannian lagrangian fibrations of almost \ka manifolds. When  one starts from mirror dual semi-flat special lagrangian tori fibrations of Calabi-Yau manifolds, the almost $2$-\ka manifolds thus obtained (actually a small deformation of them) are conjectured to be self-dual. We will prove that this is the case at least in some situations in the next  two sections.\\
The following conditions on a submersion have been already considered in the literature:
\begin{dfn}
Let $(X,\mathbf{g})$ be a Riemannian manifold, let $B$ be a smooth manifold, and let $f~:~X\to B$ be a smooth submersion.\\
1) We say that $f$ is {\em Riemannian} if there exists a (necessarily unique) Riemannian metric on $B$ such that $df$ is an isometry from $Ker(df)_p^\perp$ to $T_{f(p)}B$ for all $p\in X$.\\
2) We say that $f$ is {\em covariant constant} if it is Riemannian, and $df$ commutes with parallel transport, i.e. if $\gamma(t), t\in [0,1]$ is a path in $X$, $G_X~:T_{\gamma(0)}X\to T_{\gamma(1)}X$ is parallel transport in $X$ along $\gamma$, and $G_B~:T_{f(\gamma(0))}B\to T_{f(\gamma(1))}B$ is parallel transport in $B$ along $f(\gamma)$, then 
$G_B\left(df_{\gamma(0)})(v)\right) = df_{\gamma(1)}\left(G_X(v)\right)$ for all $v\in T_pX$.
\end{dfn}
Recall that an almost \ka manifold is a symplectic manifold together with a compatible Riemannian metric. We are now ready to state the main result of this section:
\begin{teo}
\label{specialtoska}
Let $(X_i,\omega_{X_i},\mathbf{g}_{X_i})$  be  almost \ka  manifolds of (real) dimension $2n$, for $i = 1,...,s$. Let $B$ be a smooth manifold of dimension $n$, and let $f_i:X_i\to B$  be surjections which are also lagrangian fibrations (with respect to the \ka forms).  Consider $X = X_1\times_B\cdots\times_B X_s$, with the metric $\mathbf{g}$ induced from $X_1\times\cdots\times X_s$ and with the 2-forms $(\omega_1,...,\omega_s)$, where $\omega_i$ is $\sqrt{s}$ times the pull-back of the \ka form of $X_i$, under the natural projection  $\pi_i: M\to X_i$. We then have that:\\
1) $(X,\omega_1,...,\omega_s)$ is a polysymplectic manifold.\\
2)  If all the $f_i$ are Riemannian with respect to the same metric on $B$, then $\mathbf{g}$ is compatible with the polysymplectic structure $\omega_1,...,\omega_s$. In other words,  $(M,\omega_1,...,\omega_s,\mathbf{g})$ is an almost $s$-\ka manifold.\\
3)  If moreover the $\omega_{X_i}$ are covariant constant on the respective $X_i$ and all the $f_i$ are covariant constant, then  all the $\omega_j$ are covariant constant with respect to the metric $\mathbf{g}$ on $M$.
\end{teo}
\textit{Proof}\\
1) The proof amounts to proving that the forms induce a polysymplectic structure pointwise, and that the distribution $\sum_j\omega_j^\perp$ is integrable. The first fact is an easy linear algebra observation. For the second one, let $F:X_1\times_B\cdots\times_B X_s\to B$ be the induced map, which is a fibration. Then the required integrability follows from the fact that
\[\omega_1^\perp + \omega_2^\perp = Ker(dF)\]
\\
2) Given $p\in X$, we will show that there is an orthonormal polysymplectic basis of $T_pM$. Pick an orthonormal basis $v_1,...,v_n$ of $T_{f(p)}B$, and let $z_1^j,...,z_n^j$ be a set of vectors in $Ker(d(f_j)_{p_j})^\perp$  (were $p = (p_1,...,p_s) \in X \subset X_1\times\cdots\times X_s$), such that $df_j(z_i^j) = v_i$ for all $i,j$. 
Because the $f_j$ are Riemannian, it follows that the $z_i^j$ are orthonormal (for fixed $j$). Define 
\[w_i = \frac{1}{\sqrt{s}}(z_i^1,...,z_i^s) \in T_p(X_1\times\cdots\times X_s)\] 
From their definition, it follows that the $w_i$ lie actually in $T_pM$. Moreover, $(w_l,w_m) = \delta_{lm}$. Define also 
\[w_i^j = (0,...,Jz_i^j,0,..0)~ ~(j^{th}~ place)\]
in $T_pM$. We are indicating with $J$ the almost complex structure on the various $X_i$ (or the induced one on $X_1\times\cdots\times X_s$, which is the same). The fact that the $w_i^j \in T_pM$ follows from the fact that $Jz_i^j\in Ker(d(f_j)_{p_j})$, which is a consequence of the Lagrangian condition. It is now very easy to verify that $w_1,...,w_n,w_1^1,...,w_n^s$ is an orthonormal polysymplectic basis at $p$ with respect to the polysymplectic structure $\omega_1,...,\omega_s$.\\
3) The forms $\omega_j$ are covariant constant on $X_1\times\cdots\times X_s$ (because they are by hypothesis covariant constant on their respective $X_j$'s). If all the $f_j$ are covariant constant, then parallel transport on $X$ is then the restriction of parallel transport on $X_1\times\cdots\times X_s$, and hence the $\omega_j$ are constant also on $X$.
\qed

The construction of the above theorem is natural enough to deserve a name. We will actually normalize the metric along the "horizontal" directions for a reason which will be made clear by the remark following the definition.
\begin{dfn}
\label{fibprod2k}
Let 
\[\mathbf{X_1} = (X_1,\omega_{X_1},\mathbf{g}_{X_1}),~~ \mathbf{X_2} = (X_2,\omega_{X_2},\mathbf{g}_{X_2})\]
be \ka manifolds, and let $f_i:(X_i,\me_{X_i})\to (B,\me_B)$ be smooth Riemannian surjections which are also lagrangian fibrations, with $dim(B) = n$. We then define the almost $2$-\ka manifold $\x12$ as 
\[\x12 = \left(X_1\times_B X_2,f_1^*(\omega_{X_1}),f_2^*(\omega_{X_2}),\mex12\right)\]
The metric $\mex12$ is $\alpha_{\frac{1}{\sqrt{2}}}\left(\beta_{\frac{1}{\sqrt{2}}}\left(i^*(\me_{X_1}\times\me_{X_2})\right)\right)$ where $i^*$ is the pull-back along the inclusion, $\me_{X_1}\times\me_{X_2}$ is the product metric on $X_1\times X_2$ and $\alpha_t,\beta_t$ are with respect to the almost $2$-\ka structure 
$\left(\sqrt{2}f_1^*\omega_{X_1},\sqrt{2}f_2^*\omega_{X_2},i^*(\me_{X_1}\times\me_{X_2})\right)$ given by the previous theorem
\end{dfn}
The reason we adopted the definition above is that then we have the following
\begin{rmk}
Assume that $\mathbf{X_i}$, $f_i:(X_i,\me_{X_i})\to (B,\me_B)$ are as in the previous definition. Then the naturally induced maps
\[X_1\times X_2\to X_1,~X_1\times_B X_2\to X_2,~~ X_1\times_B X_2 \to B\]
are smooth Riemannian surjections. 
\end{rmk}
One could also use the previous remark (plus the fact that the forms $\omega_i$ are pull-backs of the \ka forms on the $X_i$) to characterize $\x12$.\\
To put the condition of being Riemannian into perspective, and to make contact with Mirror Symmetry, we relate it with the semi-flatness condition of ~\cite{SYZ}, or rather with one of its consequences. We start by recalling the following standard definition:
\begin{dfn}
Let $(X,\omega,\mathbf{g},\Omega)$ be a Calabi-Yau  manifold of complex dimension $n$ (where $\omega$ is the \ka form, $\mathbf{g}$ the \ka metric and $\Omega$ the globally defined nondegenerate holomorphic $n-form$).\\
1) We say that a submanifold $L\subset X$ is {\em Special Lagrangian} if it is Lagrangian (of maximal dimension) with respect to $\omega$, and there exists a complex number of the form $e^{i\theta}$ 
such that $Im(e^{i\theta}\Omega)|_{L} = 0$. Such a $\theta$ is called the {\em phase} of the special lagrangian submanifold.\\
2) We say that a smooth map $f:X\to B$ to a smooth manifold $B$ of (real) dimension $n$ is a {\em Special Lagrangian Fibration} if $f$ is a submersion and for all $q\in B$ the submanifold $L_q = f^{-1}(q)\subset X$ is a special lagrangian submanifold of  
$(X,\omega,\mathbf{g},\Omega)$. We require also that the phase of the fibres is constant. \\
3) A special lagrangian fibration is said {\em semi-flat} if the induced metrics on the fibres are flat
\end{dfn}
\begin{pro}
Let $(X,\omega_X,\mathbf{g}_X,\Omega)_X$ be  a Calabi-Yau  manifold of  complex dimension $n$
Let $f:X\to B$ be a special lagrangian fibration with compact connected fibres, such that the metric of $X$ restricted to any fibre is flat. Then  $f$ is Riemannian.
\end{pro}
\textit{Proof}\\
In view of the description of deformations of special lagrangian manifolds of  ~\cite{ML}, it is enough to observe that harmonic forms on a flat manifold are covariant constant, and also their dual vector fields are covariant constant. As parallel transport is an isometry on any Riemannian manifold, and the complex endomorphism is also an isometry, this implies that on each fibre you have an orthonormal frame of vector fields, whose transformations under the complex involution give a complete set of first order normal deformations of the fibre itself. This clearly implies that $f$ is Riemannian.
\qed
\begin{rmk}
In the situation of the above proposition, ~\cite{ML} proves also that the map $f$ is a surjection to a smooth manifold of dimension $n$.
\end{rmk}
\begin{cor}
Let $f_1:X_1\to B$ and $f_2:X_2\to B$ be semi-flat special lagrangian fibrations of Calabi-Yau manifolds. There is then a natural structure of almost $2$-\ka manifold on $X_1\times_B X_2$. If the fibrations are covariant constant, then the forms $\omega_1$, $\omega_2$ and $\omega_D$ associated to this almost $2$-\ka structure are (Levi-Civita) covariant constant (and therefore, in particular, $(X_1\times_B X_2,\omega_1,\omega_2,\me)$ is weakly self-dual) 
\end{cor}
\begin{rmk}
One does not expect actual special lagrangian fibrations to be covariant constant, except in flat cases (arising, from example, from special lagrangian fibrations of abelian varieties)
\end{rmk} 
\begin{con}
If we start from a mirror pair $X,Y$ of Calabi-Yau manifolds, then for each point near the large complex structure limit of $X$ there is a point near the large \ka structure of $Y$ for which there are lagrangian tori fibrations of (dense open subsets of) $X$ and $Y$ over the same basis $B$ such that  there is a ("small") deformation $\mathbf{h}$ of the metric $\me$ on $X\times_B Y$ for which $(X\times_B Y,\omega_1,\omega_2, \mathbf{h})$ is self-dual, where $\omega_1$ and $\omega_2$ are $\sqrt{2}$ times the pull-backs of the \ka forms of $X$ and $Y$ respectively. 
\end{con}
\begin{rmk}
In general, we cannot expect to have fibrations of all of $X$ and $Y$ over the same $B$. We do expect  however that  $(X\times_B Y,\omega_1,\omega_2, \mathbf{h})$ admits a natural compactification (as a self-dual manifold).
\end{rmk}
Although the conjecture is not established in general, in the next two  sections we will prove it in the "limit" situation of \ak manifolds, studied for example in the papers ~\cite{KS} and ~\cite{GW}, and  then for elliptic curves (in a refined form). In these two situations the behaviour is actually simpler than the general expected one, because we have dual special lagrangian fibrations over all of $X$ and $Y$ at all mirror pairs of points from their respective moduli spaces.
\begin{rmk}
We do not expect the induced  metric on  $X\times_B Y$ to be self-dual (of even almost $2$-\ka) at finite points in the moduli spaces. The corrections necessary in order to make it self-dual could be thought  of as "quantum corrections" although, be warned, the deformed metric does not necessarily induce one either on $X$ or on $Y$.\\
The (conjectural) recipe to find the self-dual manifolds associated to mirror pairs is then the following: find "mirror dual" lagrangian fibrations on big enough open subsets $X$ (resp. $Y$) of the manifolds over the same basis $B$  by going near large complex (resp. \ka) structure points of the moduli spaces. Put a self-dual metric on $X\times_B Y$, and then compactify what you obtained in the category of self-dual manifolds.
\end{rmk}
\section{\ak manifolds}
\label{sec:maandell}
In this seciton we will show how one can build self-dual manifolds starting from \ak manifolds. We first need to recall a definition and two lemmas from ~\cite{KS}:
\begin{dfn}[\cite{KS}, Definition 2, page 17]
An \ak manifold is a triple $(Y,g,\nabla)$ where $(Y,g)$ is a smooth Riemannian manifold with the metric $g$, and $\nabla$ is a flat connection on $TY$ such that:\\
a) $\nabla$ defines an affine structure on $Y$\\
b) Locally in affine coordinates $(x_1,...,x_n)$ the matrix $(g_{ij})$ of $g$ is given by $g_{ij} = \partial^2K/\partial x_i\partial x_j$ for some smooth real-valued function $K$ on $Y$.\\
If moreover one has that\\ 
c) The \ma equation $det(g_{ij}) = const$ is satisfied\\
then $(Y,g,\nabla)$ is called an {\em \ma } manifold.
\end{dfn}
\begin{lem}[\cite{KS}, Proposition 2, section 3.2]
For a given \ak manifold $(Y,g_Y,\nabla_Y)$ there is a canonically defined dual \ak manifold $(Y^\vee,g_{Y^\vee},\nabla_{Y^\vee})$ such that $(Y,g_Y)$ is identified with $(Y^\vee,g_{Y^\vee})$ as Riemannian manifolds, and the local system $(T_{Y^\vee},\nabla_{Y^\vee})$ is naturally isomorphic to the local system dual to $(TY,\nabla_Y)$.
\end{lem}
\begin{lem}[\cite{KS}, Corollary 1, section 3.2] 
If $\nabla_Y$ defines an integral affine structure on $Y$, then $\nabla_{Y^\vee}$ defines an integral affine structure on $Y^\vee$. As the dual covariant lattice one takes the lattice $(TY^\znum)^\vee$ dual to $TY^\znum$ with respect to the metric $g_Y$.
\end{lem}
\begin{teo}
Let $(Y,g_Y,\nabla_Y)$ be a \ak manifold, such that $\nabla_Y$ defines an integral affine structure. Let $(Y^\vee,g_Y^\vee,\nabla_{Y^\vee}) = (Y,g_Y,\nabla_Y^\vee)$ be its dual \ak manifold. We then have that there is a canonically induced almost $2$-\ka structure $(\omega_1,\omega_2,\mathbf{h})$ on 
\[X_Y = \left(TY/TY^\zint\right)\times_Y \left(TY/(TY^\zint)^\vee\right)\]
and with this structure $X_Y$ is self-dual
\end{teo}
\textit{Proof}
To build the almost $2$\ka structure we first put a Riemannian metric on $TY$, using the flat connections given by the affine structure to select the orthogonal complements to the fibres of the projection to $Y$.  The metric $\me = \me_Y$ on  $Y$ then induces via the projection the metric on these horizontal distributions, and by translationn that on the fibre directions. By construction we get that the fibration $TY\to Y$ is Riemannian. The \ka form is the pull-back of the standard symplectic form of $T^*_Y$ to $TY$ via the map induced by the metric on $Y$. If we choose coordinates $x_1,...,x_n$ on (an open set inside) $Y$, we have induced coordinates $(x_1,...,x_n,y_1 = dx_1,...,y_n = dx_n)$ on $TY$, and in these coordinates the symplectic form is
\[\omega_{TY} = \sum_{i,k} \me_{ik}(x_1,...,x_n)dx_i\wedge dy_k\]
It is immediate to verify that this form is compatible with the metric, and defines an almost complex structure which is integrable (this \ka structure  can be identified with that of ~\cite{KS}, paragraph 5.2, once we identify $TY$ with $T^*_Y$ via the metric $\me$).\\ 
The projection to $Y$ is Riemannian and  lagrangian. This is all that is needed in  Theorem ~\ref{specialtoska} and Definition ~\ref{fibprod2k}, and we get therefore an almost $2$-\ka structure on $TY\times_Y TY$, with $\omega_i = \pi_i^*\omega_{TY}$, where $\pi_i$ are the projections on the two factors, and with metric $\mathbf{h}$ induced from $\me\times\me$ on  $TY\times TY$ in the way described in Definition ~\ref{fibprod2k}. Let $\pi_Y:TY\times_B TY \to Y$ be the canonical projection. We then have that $TY\times_B TY$ is a vector bundle on $Y$, and the metric $\mathbf{h}$ and the forms $\omega_1,\omega_2$ are invariant with respect to translations by covariant constant sections. As the integral lattices are generated by covariant constant sections, it follows that we get an almost $2$-\ka structure also on the quotient 
\[X_Y = TY\times_Y TY/\left(\pi_1^{-1}(TY^\zint)\times \pi_2^{-1}({TY^\zint}^\vee)\right)\]
We continue to call $\pi_Y$ the projection from $X_Y$ to $Y$, which is Riemannian also with respect to $\mathbf{h}$, and induces on $Y$ the metric $\me$.\\
We now choose integral affine coordinates $x_1,...,x_n$ on (an open set inside ) $Y$, and indicate with $y^1_1,...,y^1_n$ (resp. $y^2_1,...,y^2_n$) the induced coordinates on the first copy (resp. the second copy) of $TY$. This determines coordinates 
\[x_1,...,x_n,y^1_y,...,y^1_n,y^2_1,...,y^2_n\]
on $TY\times_Y TY$, which can also be used locally on $X_Y$. With these coordinates,
\[\omega_j = \sum_{ik}\me_{ik}(x_1,...,x_n)dx_i\wedge dy^j_k\]
The affine coordinates $z_1,...,z_n$ on $Y$ dual to $x_1,...,x_n$ satisfy (by definition) 
\[\frac{\partial z_i(x_1,...,x_n)}{\partial x_k} = \me_{ik}(x_1,...,x_n)\]
and therefore the coordinates $w_1,...,w_n$ in the fibre directions associated to $z_1,...,z_n$ satisfy the relation $w_k = \sum_i\me_{ik}y_i$.
It follows that the leaves of the horizontal distribution associated to the connection dual to $\nabla_Y$ are described (locally) by $\left\{(x_1,...,x_n,\sum_i\me^{i1}w_i,...,\sum_i\me^{in}w_i)\right\}$, for numbers $w_1,...,w_n$. The dual horizontal distribution at the point $(x_1,...,y_n)$ of $TY$ is therefore generated by the vectors
\[\frac{\partial}{\partial x_i} + \sum_{lkm}y_m\me_{mk}\frac{\partial\me^{lk}}{\partial x_i}\frac{\partial}{\partial y_l}\]
The metric $\mathbf{h}$  makes the vectors 
$v_i = (\frac{\partial}{\partial x_i} + \sum_{lkm}y_m\me_{mk}\frac{\partial\me^{lk}}{\partial x_i}\frac{\partial}{\partial y_l^2})$ orthogonal to the fibre directions of both the projections $\pi_1$ and $\pi_2$. Moreover,
\[\omega_1(\frac{\partial}{\partial y^1_i},v_k) = -\me_{ik} = \omega_2(\frac{\partial}{\partial y^2_i},v_k) \]
and therefore the map from $\omega_2^\perp$ to $\omega_1^\perp$ induced by the almost $2$-\ka structure is simply induced by the correspondence $\frac{\partial}{\partial y^1_i}\to \frac{\partial}{\partial y^2_i}$. The form $\omega_D$ is simply the differential form associated to the composition of this map with the natural map in the dual of $\omega_1^0$ built with  the metric, and therefore it will be of the form $\sum_idy^1_i\wedge \alpha_i$, where $\alpha_i$ is the $1$-form annichilating  $\omega_2^\perp + V$ and dual to $\frac{\partial}{\partial y^2_i}$ inside $\omega_1^0$, where $V$ is the orthogonal complement to the subspace $\omega_1^\perp + \omega_2^\perp$. For the choice
\[\alpha_i = \sum_j\me_{ij}dy^2_j - \sum_{lmkj}y^2_m\me_{mk}\me_{il}\frac{\partial \me^{lk}}{\partial x_j}dx_j\]
\[\alpha_i(v_h) = \sum_{lkm}\me_{il}y^2_m\me_{mk}\frac{\partial \me^{lk}}{\partial x_h} - \sum_{lmk}y^2_m\me_{mk}\me_{il}\frac{\partial \me^{lk}}{\partial x_h} = 0, ~\alpha_i(\frac{\partial}{\partial y^2_j}) = \me_{ij}\]
and therefore 
\[\omega_D = \sum_i dy^1_i\wedge \alpha_i = \sum_{ij}\me_{ij}dy^1_i\wedge dy^2_j - \sum_{ijlmk}y^2_m\me_{mk}\me_{il}\frac{\partial \me^{lk}}{\partial x_j}dy^1_i\wedge dx_j\]
To show that $X_Y$ is weakly self-dual, we must prove that $d\omega_D = 0$. 
\[\begin{array}{l}
d\omega_D = \sum_{ijk}\frac{\partial\me_{ij}}{\partial x_k}dy^1_i\wedge dy^2_j\wedge dx_k - \sum_{ijklm}\me_{mk}\me_{il}\frac{\partial \me^{lk}}{\partial x_j} dy^2_m\wedge dy^1_i\wedge dx_j -\\
\sum_{ijklmn}y^2_m\left(
\frac{\partial \me_{mk}}{\partial x_n}\me_{il}\frac{\partial \me^{lk}}{\partial x_j} + 
\me_{mk}\frac{\partial \me_{il}}{\partial x_n}\frac{\partial \me^{lk}}{\partial x_j} +
\me_{mk}\me_{il}\frac{\partial^2 \me^{lk}}{\partial x_j\partial x_n} 
\right)dy^1_i\wedge dx_j\wedge dx_n 
\end{array}\]
By changing name to the summation indices the vanishing of the first line is equivalent to that of 
\[\sum_{ij}\left(\frac{\partial\me_{ij}}{\partial x_k} + \sum_{lm}\me_{jm}\me_{il}\frac{\partial \me^{lm}}{\partial x_k}\right)dy^1_i\wedge dy^2_j \]
Multiplying with $\me^{jr}$ and summing over $j$ we get that the vanishing of the first line is equivalent to the vanishing of
$\sum_j\me^{jr}\frac{\partial\me_{ij}}{\partial x_k} + \sum_{l}\me_{il}\frac{\partial \me^{lr}}{\partial x_k} = \frac{\partial}{\partial x_k}\left(\sum_j \me^{jr}\me_{ij}\right) = \frac{\partial}{\partial x_k}\delta_{ir}$
and this is clearly zero.\\
Coming to the second line, it is clear that we need the vanishing of
\[\sum_{jn}\sum_{kl}\left(
\frac{\partial \me_{mk}}{\partial x_n}\me_{il}\frac{\partial \me^{lk}}{\partial x_j} + 
\me_{mk}\frac{\partial \me_{il}}{\partial x_n}\frac{\partial \me^{lk}}{\partial x_j} +
\me_{mk}\me_{il}\frac{\partial^2 \me^{lk}}{\partial x_j\partial x_n} 
\right)dx_j\wedge dx_n\]
We are therefore reduced to proving the symmetry in $j,n$ of
\[\sum_{kl}\left(
\frac{\partial \me_{mk}}{\partial x_n}\me_{il}\frac{\partial \me^{lk}}{\partial x_j} + 
\me_{mk}\frac{\partial \me_{il}}{\partial x_n}\frac{\partial \me^{lk}}{\partial x_j}
\right)\]
By multiplying by $\me^{ir}\me^{ms}$ and summing over $i$ and $m$ we reduce this to
\[\sum_{lkim}\left(\me^{ms}\me^{ir}\frac{\partial \me_{mk}}{\partial x_n}\me_{il}\frac{\partial \me^{lk}}{\partial x_j} + \me^{ms}\me^{ir}\me_{mk}\frac{\partial \me_{il}}{\partial x_n}\frac{\partial \me^{lk}}{\partial x_j}\right) = 
\sum_{km}\me^{ms}\frac{\partial \me_{mk}}{\partial x_n}\frac{\partial \me^{rk}}{\partial x_j} 
+ \]
\[\sum_{li}\me^{ir}\frac{\partial \me_{il}}{\partial x_n}\frac{\partial \me^{ls}}{\partial x_j} = 
- \sum_{km}\me_{mk}\frac{\partial \me^{ms}}{\partial x_n}\frac{\partial \me^{rk}}{\partial x_j} - \sum_{li}\me_{il}\frac{\partial \me^{rl}}{\partial x_n}\frac{\partial \me^{ls}}{\partial x_j}\]
and the last expression is clearly symmetric in $j,n$.\\
We have therefore that $(X_Y,\omega_1,\omega_2,\mathbf{h})$ is weakly sefl-dual. To prove that it is self-dual it remains to be shown that the fibres of the projection to $Y$ have all Riemannian volume one. This however is clear, as they are all of the form $T\times T^\vee$, where $T$ is a torus and $T^\vee$ is its dual with respect to the metric, and it is a general fact that in this case $vol(T\times T^\vee) = vol(T)vol(T^\vee) = 1$.
\qed
\section{Elliptic curves}
\label{sec:ellcur}
Before going into the characterization of mirror symmetry for elliptic curves in terms of self-dual manifolds, we need to define it. In this case there are no ambiguities, and all is clear and settled by now. First we recall some terminology from ~\cite{D1} (or equivalently from ~\cite{PZ} or ~\cite{D2}).
\begin{dfn}[See for example ~\cite{D1}, pages 152-153]
Let $(\tau,t)\in \mathbb{H}\times\mathbb{H}$, where $\mathbb{H}$ is the upper half plane inside $\cnum$. We associate to $(\tau,t)$ the complex manifold  $E_\tau = \cnum/\znum\tau\oplus\znum$  with the  complexified \ka form
\[\omega^{\tau,t} = -\frac{t}{2Im(\tau)}dz\wedge d\bar{z}\]
Such a pair $(E_\tau,\omega^{\tau,t})$ is also indicated with $\mathbf{E}_{\tau,t}$.
\end{dfn}
Notice that our complexified \ka class must be multiplied by $2\pi$ to recover that of \cite{D1}. The imaginary part of $\omega^{\tau,t}$ (multiplied by $2\pi$) is what is usually called the $B$-field, while the real part is a \ka form on $E_\tau$.
\begin{dfn} The elliptic curve  with a complexified \ka class $E_{\tau,t}$ is mirror dual to the  elliptic curve with a complexified \ka class $E_{t,\tau}$
\end{dfn}
For a justification of the above definition, see for example ~\cite{PZ} or ~\cite{D1} and the references therein. We will not get into this justification here.
\begin{rmk}
The natural projection map $E_\tau\to B = S^1$ induced by the projection $\cnum\to i\rnum$ is (special) lagrangian and Riemannian (with respect to the flat metric $\frac{t_2}{\tau_2}(dxdx + dydy)$ on $\cnum$). With the induced metric the basis has length $\sqrt{t_2\tau_2}$. 
\end{rmk}
\begin{dfn}
Let $(\tau,t) = (\tau_1 + i\tau_2,t_1 + it_2) \in \mathbb{H}\times\mathbb{H}$, where $\mathbb{H}$ is the upper half plane inside $\cnum$. We associate to $(\tau,t)$ the almost $2$-\ka manifold 
\[\mathbf{X}_{\tau.t}= (X_{\tau,t},\omega_1^{\tau,t},\omega_2^{\tau,t},\me^{\tau,t})\]
defined as $X_{\tau,t} = E_\tau\times_B E_t$, with the almost $2$-\ka structure induced as in Definition ~\ref{fibprod2k} by the \ka structures $-i\frac{t_2}{2\tau_2}dz\wedge d\bar{z}$ and $-i\frac{\tau_2}{2t_2}dz\wedge d\bar{z}$  on $E_\tau$ and $E_t$ respectively.
\end{dfn}
\begin{lem}
\label{xtautselfdual}
For any choice of $(\tau,t)\in \mathbb{H}\times\mathbb{H}$, the almost $2$-\ka manifold $\mathbf{X}_{\tau,t}$ is self-dual
\end{lem}
\textit{Proof}
The forms $\omega_1$ and $\omega_2$ are covariant constant, therefore the manifold is automatically weakly self-dual. To check that the leaves of $\omega_1^\perp+\omega_2^\perp$ have Riemannian volume one, observe that the leaves of $E_\tau\to B$ are of the form $\rnum/\znum$ with metric $\frac{t_2}{\tau_2}dx$, while those of  $E_t\to B$ are of the form $\rnum/\znum$ with metric $\frac{\tau_2}{t_2}dx$. This proves that the volume of the product is one, as desired.
\qed
\begin{rmk}
Let $\mathbf{X} = (X,\omega_1,\omega_2,\me) \cong \mathbf{X}_{\tau,t}$ as almost $2$-\ka manifold. We then have that the (a priori non-commutative) quotients 
\[E_1 = X/\omega_1^\perp,~~E_2 = X/\omega_2^\perp,~~B = X/\omega_1^\perp+\omega_2^\perp\]
are smooth manifolds, and the natural projection maps $\pi_i:X\to E_i$ are smooth Riemannian. Moreover, $(E_1,{\pi_1}_*\omega_1,{\pi_1}_*\me)$ and  $(E_2,{\pi_2}_*\omega_1,{\pi_2}_*\me)$ are elliptic curves both fibred (with lagrangian Riemannian maps) onto $B$. $\mathbf{X}$ can be recovered as $E_1\times_B E_2$ (with the induced almost $2$-\ka structure).
\end{rmk}
\begin{lem}
The fibration $X_{\tau,t}/\omega_2^\perp\to X_{\tau,t}/\omega_1^\perp+\omega_2^\perp$ is a principal fibration with group $S^1$ and monodromy $t_2\in S^1$ around the generator of $\pi_1$ of the basis. Similarly, the fibration $X_{\tau,t}/\omega_1^\perp\to X_{\tau,t}/\omega_1^\perp+\omega_2^\perp$ is a principal fibration with group $S^1$ and monodromy $t_1\in S^1$ around the generator of $\pi_1$ of the basis. 
\end{lem}
\textit{Proof}
We are simply considering the fibration $\cnum/\znum+t\znum\to i\rnum/it_2\znum$ induced by the fibration $\cnum\to i\rnum$, with a multiple of the flat metric. The statement is then clear. The second statemet is proved  in the same way.
\qed
\begin{lem}
Let $\mathbf{X} = \mathbf{X}_{\tau,t}$ be the self-dual manifold associated to the pair $(\tau,t)$. We then have that:\\
1) The lenght of the manifold $X/\omega_1^\perp+\omega_2^\perp$ (with the induced metric) is $\sqrt{t_2\tau_2}$\\
2) The lenght of the fibre of the fibration $X\to X/\omega_1^\perp$ is $\sqrt{\frac{t_2}{\tau_2}}$\\
2) The lenght of the fibre of the fibration $X\to X/\omega_2^\perp$ is $\sqrt{\frac{\tau_2}{t_2}}$
\end{lem}
\textit{Proof}
All the statemenst are easy calculations. We omit the deails.
\qed

The lemmas above imply the following
\begin{teo}
\label{sddetcur}
The self-dual manifold $X_{\tau,t}$ determines (by a well defined procedure) the elliptic curves with a complexified \ka class $\mathbf{E}_{\tau,t}$ and $\mathbf{E}_{t,\tau}$
\end{teo}
The result above should not come as a surprise, as the self-dual manifold determines the elliptic curves toghether with a metric and a special lagrangian fibration.
\begin{dfn}
We will indicate with 
\[\mathbf{X}_{\tau,t}\to_1 \mathbf{E}_{\tau,t},~~\mathbf{X}_{\tau,t}\to_2 \mathbf{E}_{t,\tau}\]
the content of the previous theorem
\end{dfn}
This definition allows us to state more precisely what has been proved:
\begin{teo}
Let $\mathbf{E}_1,\mathbf{E}_2$ be elliplic curves together with complexified \ka classes. Then the following are equivalent:\\
1) $\mathbf{E}_1$ and $\mathbf{E}_2$ form a  mirror pair\\
2) There is a self-dual manifold of the form $\mathbf{X}_{\tau,t}$ such that 
\[\mathbf{X}\to_1 \mathbf{E}_1, ~ ~ \mathbf{X}\to_2 \mathbf{E}_2\]
\end{teo}
From our point of view, the situation of elliptic curves (and very likely of abelian varieties in general) is a degenerate one, in which the description in terms of self-dual manifolds and that in terms of $B$-fields are equivalent. In general, we expect that the description in terms of $B$-fields and special lagrangian fibrations holds only "in the limit", while self-dual manifolds exist also at finite points, and converge to the limit situation near the boundary of the moduli space.\\
The following remark might be useful to recover another part of the classical terminology.
\begin{rmk}
\[(X_{\tau,t},\me^{\tau,t})\to (E_{\tau},\me_{E_{\tau,t}})\]
in the sense of Gromov-Hausdorff, for $\frac{Im(\tau)}{Im(t)}\to +\infty$,
\end{rmk}
For reasons of space (and time) we do not analyze in detail what happens for elliptic fibrations of $K3$ surfaces, which (together with abelian varieties) would be the next natural step to take. However, to give the reader something to think about, we formulate a simple conjecture which relates the above constructions to those of ~\cite{GW}.\\
Let $X\to B$ and $\hat{X}\to B$ be mirror dual  semi-flat special lagrangian fibrations of a general mirror pair of $K3$ surfaces over the same basis $B$ (take your favourite definition for what that is). This fibration will have for a general $K3$ exactly $24$ singular fibres. Call $B^0$ the complement of the singular set inside $B$. As we have proven in Theorem ~\ref{specialtoska}, we then have an almost $2$-\ka structure on $X\times_{B^0}\hat{X}$. 
\begin{con}
In the situation described above, and with the induced almost $2$-\ka structure,  $X\times_{B^0}\hat{X}$ is self-dual. Moreover, for all $b\in B$ the class of the form $\omega_D$ restricted to the fibre over $b$ is a constant (independent of $b$) multiple of the first Chern class of the Poincar\'{e} bundle of that fibre (once you rotate the complex structure of $X$ and $\hat{X}$ to make the fibrations to $B$ holomorphic). The self-dual manifold $X\times_{B^0}\hat{X}$ admits a natural compatictification to a smooth compact self-dual manifold.
\end{con}
This conjecture should shed some light on the nature of the form $\omega_D$, and on why we expect it to be closed in situations arising from mirror pairs.
\section{Polysymplectic and $2$-\ka manifolds}
\label{sec:polyandska}
In this section we describe some of what happens if we weaken (in the polysymplectic case) or strengthen (in the $2$-\ka case) the condition of self-duality. These two notions were introduced in ~\cite{G}, where the reader can go to find a more detailed study. As the following will be general considerations, we will not need to stick to the case $s=2$. Remenber however that we defined self-dual manifolds only for $s=2$ (although it would be easy to generalize the definition to general $s>1$).
\begin{teo}[Polysymplectic normal form]
\label{teo-multisympnorm}
Let $(X,\omega_1,...,\omega_s)$ be a smooth polysymplectic manifold and $p\in X$. Assume given elements $\phi_1,...,\phi_n,\psi^1_1,...,\psi^n_s$ of $T^*_pX$ such that for all $j=1,..,s$ one has $(\omega_j)_p=\sum_i\phi_i\wedge\psi^i_j$.Then there are a neighborhood $\mathcal{U}\subset X$ of $p\in X$, a neighborhood $\mathcal{V}\subset \mathbf{R}^{dim(X)}$ of $0\in \mathbf{R}^{dim(X)}$ and an isomorphism of polysymplectic manifolds
\[\phi~:~(\mathcal{U},\omega_1,...,\omega_s)~\rightarrow~
\left(\mathcal{V},\sum_i dx_i\wedge dy^1_i,...,\sum_i dx_i\wedge dy^s_i\right)\]
where we indicated the coordinates on $\mathbf{R}^{dim(X)}$ with $x_1,...,x_n,y^1_1,...,y^s_n$. With this notation, one can also assume $(dx_i)_p=\phi$, $(dy^i_j)_p=\psi^i_j$.
\end{teo} 
{\it Proof}\\
If $V$ is a vector space, given an element of $\alpha\in\bigwedge^*(V)$ we indicate with $C(\alpha)$ the smallest subspace $W\subset V$ such that $\alpha\in \bigwedge^*W$. Similarly, for a differential form $\alpha$ we define $C(\alpha)$ to be the smallest distribution of subspaces $D\subset\Omega^1$ such that $\alpha\in \bigwedge^*D$. 
A priori, the $C(\omega_j)$ are only "generalized Pfaffian systems", as defined for example in ~\cite[Page 382]{LM}. 
From Darboux's Reduction Theorem, in the form stated for example in ~\cite[Bryant, Page 103]{FU}, we see that $C(\omega_j)$ is a vector bundle (of rank $2n$) for any $j=1,...,s$, with local coframes given by closed $1$-forms. Define $C^X=\bigcap_j C(\omega_j)$. Then $C^X = \left(\sum_j\omega_j^\perp\right)^\perp$ is a constant rank distribution of subspaces of $T^*X$, and by the definition of a polysymplectic structure and Frobenius it is locally generated by closed forms. From the constant rank property, we may assume that there are (locally) $n$ functions $x_1,...,x_n$ such that $dx_1,...,dx_n$ are independent, and for all $q$ in the open set considered
\[< (dx_1)_p,\dots,(dx_1)_p > = C^X\]
By acting if necessary with a constant transformation matrix we can assume that $\forall i~ (dx_i)_p=\phi_i$.\\
Fix now an index $j\in\{1,...,s\}$. From Darboux's reduction theorem, we can find coordinates $z_1,...,z_d$ such that $\omega_j$ is expressed only in terms of $z_{d-2n+1},...,z_d$, and such that $\frac{\partial}{\partial z_k}$ is in $C(\omega_j)^\perp$ for $k=1,..,d-2n$ (and therefore one has also $<dz_{d-2n+1},...,dz_d>=C(\omega_j)$). From their definition, it follows that $\frac{\partial x_i}{\partial z_k}=0$ for all $i$, and for $k=1,..,d-2n$. Therefore, we can apply the theorem of Carath\'{e}odory-Jacobi-Lie (see ~\cite[Page 136]{LM}) to conclude that there are functions $y_i^j$ (depending only on the $z_{d-2n+1},...,z_d$)  such that $dy^j_i\in C(\omega_j)$ and $\omega_j~=~
\sum_i dx_i\wedge dy^j_i$. Because $<dx_1,...,dx_n,dy^1_j,...,dy^n_j>=C(\omega_j)$, by an invertible linear transformation inside $C(\omega_j)$ (with constant coefficients) leaving all the $dx_i$ fixed  we can also assume that $dy^i_j=\psi^i_j$. After repeating the procedure for all $j$, we end up with functions $x_1,..,x_n,y^1_1,...,y^s_n$ near $p\in X$.
The $x_1,...,x_n,y^1_1,...,y^s_n$ form a system of coordinates because the $dx_1,...,dx_n,dy^1_1,...,dy^s_n$ are independent forms.
\qed

One could try to give a more conceptual proof, similar to Moser's proof of the theorem of Darboux for symplectic manifolds, using the tecniques of ~\cite{G}. This however would have taken us too far away from the theme of the present work.
\begin{cor}
Let $M$ be a smooth manifold, and $\omega_1,...,\omega_2$ be smooth $2$-forms on it. The following are then equivalent:\\
1) $(M,\omega_1,...,\omega_s)$ is a polysymplectic manifold.\\
2) For all $p\in M$ there are coordinates $x_1,..,x_n,y^1_1,...,y^s_n$ near $p$ such that 
\[\forall j~ ~\omega_j = \sum_i dx_i\wedge dy^j_i\]
\end{cor}
This characterization of polysymplectic manifolds makes clear why we consider them a natural generalization of symplectic ones.
\begin{teo}
Let $(M,\omega_1,...,\omega_s)$ be a polysymplectic manifold. The space of Riemannian metrics on $M$ compatible with the polysymplectic structure is non-empty and contractible.
\end{teo}
For the purposes of this proof, we give the following definition. 
\begin{dfn}
Let $(V,\omega_1,...,\omega_s)$ be a vector space with a non-degenerate polysymplectic structure, $s > 1$. A Riemannian metric $\mathbf{g}$ on $V$ is {\em block-compatible} with the polysymplectic structure it there exists a polysymplectic basis $e_1,...,e_n,$ $f^1_1,...,f^s_n$ such that for all $i,m,j,k$ (with $j\not= k$)
\[\mathbf{g}(e_i,f^j_m) = \mathbf{g}(f^k_i,f^j_m) = 0\]
\end{dfn}
\begin{lem}
Let $(V,\omega_1,...,\omega_s)$ be a vector space with a  polysymplectic structure, $s > 1$, and let $\mathbf{g}_1$ and $\mathbf{g}_2$ be two Riemannian metrics on $V$ block-compatible with the polysymplectic structure, and such that their restrictions to the span of the spaces $\omega_j^\perp$  coincide. If $t\in [0,1]$, then the Riemannian metric $t\mathbf{g}_1 + (1-t)\mathbf{g}_2$ is also block-compatible with the polysymplectic structure.
\end{lem}
\textit{Proof}\\
In view of the block-compatibility of the two metrics with the polysymplectic structure
there are vectors $d_1,...,d_n$, $f_1,...,f_n$, $f^1_1,...,f^s_n$ such that $e_1,...,e_n,f^1_1,...,f^s_n$ and $d_1,...,d_n,f^1_1,...,f^s_n$ are polysymplectic bases, and for $j \not= k$
\[\mathbf{g}_1(e_i,f^j_m) = \mathbf{g}_1(f^k_i,f^j_m) = 0, \mathbf{g}_2(d_i,f^j_m) = \mathbf{g}_2(f^k_i,f^j_m) = 0\]
Moreover,  we can take for all $j$ bases $h^j_1,...,h^j_n$ of the span of $f^j_1,...,f^j_n$, orthonormal with respect to $\mathbf{g}_1$ (and therefore also with respect to $\mathbf{g}_2$). We do not require such bases $h^j_1,...,h^j_n$ to be part of a polysymplectic basis. Such a basis exists because of the hypothesis on the behaviour of the two metrics on the span of $f^1_1,...,f^s_n$.
Observe first that if we define the vectors 
\[f_i(t) = e_i + \sum_{k,m}(t-1)\mathbf{g}_2(e_i,h^k_m)h^k_m,\]
then for all $i,j,m$
\[(t\mathbf{g}_1+ (1-t)\mathbf{g}_2)(f_i(t),h^j_m)  = 0\]
We now observe that there must be $\eta^m_{ik}$ such that $d_i = e_i + \sum_{km}\eta^m_{ik}h^k_m$. From the fact that 
$\mathbf{g}_2(d_i,h^k_m) = 0$, 
we deduce that $\eta^m_{ik} = - \mathbf{g}_2(e_i,h^k_m)$.
This shows that $f_i(t) = te_i + (1-t)d_i$ for all $i$, or in other words $f_i(t) = e_i + (t-1)\sum_{km}\eta^m_{ik}h^k_m$, from which it is easy to deduce that  $f_1(t),...,f_n(t),f^1_1,...,f^s_n$ is a polysymplectic basis for all $t$. This polysymplectic basis shows that $t\mathbf{g}_1+ (1-t)\mathbf{g}_2$ is block-compatible with the polysymplectic structure.
\qed
\begin{lem}
\label{existblockcomp}
Let $(M,\omega_1,...,\omega_s)$ be a (non-degenerate) polysymplectic manifold. There exists a Riemannian metric on $M$ block-compatible point by point with the polysymplectic structure.
\end{lem}
{\em Proof}\\
Pick a covering of $M$ by polysymplectic coordinate sets $\mathcal{U}_\alpha$, and a partition of unity $\{f_\alpha\}$ subordinated to the covering.\\
Observe first that if $\mathbf{g}_1$ and $\mathbf{g}_2$ are two Riemannian metrics on $M$ such that for all points $p\in M$ and for any polysymplectic basis $e_1,...,e_n,f^1_1,...,f^s_n$ of $T_pM$, for $j\not = k$, 
$\mathbf{g}_1(f^k_i,f^j_m) = 0 = \mathbf{g}_2(f^k_i,f^j_m) = 0,$
then also $t\mathbf{g}_1+ (1-t)\mathbf{g}_2$ has this property. Therefore,  by using the polysymplectic coordinates on the sets 
$\mathcal{U}_\alpha$, and the partition of unity to sum, we can easily define a Riemannian metric $\mathbf{g}$ on all of $M$ which has the property above at all points $p\in M$.
Define now a family $\mathbf{g}_\alpha$ of block-compatible metrics on any fixed open set  $\mathcal{U}_\alpha$, with the property that $\mathbf{g}_\alpha$ coincides with the fixed $\mathbf{g}$ on the span of $f^1_1,...,f^s_n$ for some, and therefore any, polysymplectic basis. Using the partition of unity, and the previous lemma, we see that we can sum all these metrics to provide a globally defined block-compatible Riemannian metric.
\qed

\begin{lem}
\label{datacompatib}
Let $(M,\omega_1,...,\omega_s)$ be a (non-degenerate) polysymplectic manifold. There is then a one to one correspondence between the following data:\\
1) A Riemannian metric on $M$, compatible with the polysymplectic structure.\\
2) A positive definite non degenerate symmetric bilinear form $\mathbf{g}^1$ on $\bigcap_{j>1}\omega_j^\perp$, plus a constant rank distribution of subspaces $W$ of $TM$, such that at each point $p\in M$ and for some polysymplectic basis $e_1,...,e_n,f^1_1,...,f^s_n$ of $T_pM$, $\mathbf{g}^1|_{T_pM}$ is supported on the span of $f^1_1,...,f^1_n$,  and $W_p = <e_1,...,e_n>$.\\
In the direction from $1)$ to $2)$ the correspondence sends a metric $\mathbf{g}$ to the bilinear form $\mathbf{g}^1$ and the subspace $W$ defined for any $p$ and any polysymplectic basis  $e_1,...,e_n,f^1_1,...,f^s_n$ of $T_pM$ as $\mathbf{g}^1|_{T_pM} = \mathbf{g}|_{<f^1_1,...,f^1_n>}$ and $W_p = <f^1_1,...,f^s_n>^{\perp_\mathbf{g}}$ respectively
\end{lem}
\textit{Proof}\\
In the direction from $1)$ to $2)$, to check that the correspondence is well defined it is enough to observe that $W_p = <e_1,...,e_n>$ for any orthonormal polysymplectic basis  $e_1,...,e_n,f^1_1,...,f^s_n$.
In the direction from $2)$ to $1)$, to define $\mathbf{g}|_{T_pM}$ choose any polysymplectic basis $e_1,...,e_n,f^1_1,...,f^s_n$ such that $W_p = <e_1,...,e_n>$, and 
$f^1_1,...,f^1_n$ is $\mathbf{g}^1|_{T_pM}$-orthonormal. Then declare any such basis to be $\mathbf{g}$-orthonormal. To check that this definition is correct, suppose given any other polysymplectic basis with the same property. Then it is immediate to check, using the observation that if a matrix is orthogonal also the transpose of its inverse is so (and actually coincides with it), that the transition matrix from one basis to the other is orthogonal, and therefore $\mathbf{g}$ is well defined. By construction, the metric $\mathbf{g}$ is Riemannian, and compatible with the polysymplectic structure point by point. The verification that the metric defined varies smoothly as $p$ varies in $M$ is straightforward, and left to the reader.
Both the correspondences thus defined are one to one and onto, as they are one the inverse of the other.
\qed

\textit{Proof of the theorem}\\
Pick any globally defined block-compatible Riemannian metric $\mathbf{g}_0$ on $M$, which exists from Lemma ~\ref{existblockcomp}. At any given point $p\in M$, pick any polysymplectic basis $e_1,...,e_n,f^1_1,...,f^s_n$, and consider the bilinear form $\mathbf{g}^1|_{T_pM} = \mathbf{g}|_{<f^1_1,...,f^1_n>}$ and the subspace $W_p = <f^1_1,...,f^s_n>^{\perp_\mathbf{g}}$. The bilinear form  $\mathbf{g}^1$ and the distribution of subspaces $W$ thus defined determine uniquely a Riemannian metric compatible with the polysymplectic structure, in view of Lemma \ref{datacompatib}.\\
To see that the space of compatible metrics is contractible, pick any metric $\mathbf{g}_0$ in it. Using Lemma \ref{datacompatib}, it is easy to see that there is a canonical way to interpolate between 
$\mathbf{g}_0$ and any other metric $\mathbf{g}$ compatible with the polysymplectic structure, and that this interpolation procedure provides a retraction of the space of compatible metrics to its point $\mathbf{g}_0$.
\qed

\begin{dfn}[\cite{G}, Definition 7.2 page 35]
A smooth manifold $M$ of dimension $n(s+1)$ together with a  Riemannian metric $\mathbf{g}$ and $2$-forms  $\omega_{1},...,\omega_{s}$ is {\em s-\ka} if the data satisfies the  following property: For each point of $M$ there exist an open neighborhood $\mathcal{U}$  of $p$ and a system of coordinates  $x_{i},y^{j}_{i}$,$i~=~1,...,n$, $j~=~1,...,s$ on $\mathcal{U}$ such that:\\
1) $\forall j ~ ~\omega_{j}~=~\sum_{i}dx_{i}\wedge dy^{j}_{i}$,\\
2) $\mathbf{g}_{(\mathbf{x},\mathbf{y})}~=~\sum_{i}dx_{i}\otimes  dx_{i}~+~\sum_{i,j}dy^{j}_{i}\otimes d 
y^{j}_{i}~+~\mathbf{O}(2)$.\\
Any such system of coordinates is called {\em standard}(s-\ka).
\end{dfn}
\begin{teo}
Let $(X,\omega_1,...\omega_s)$ be a polysymplectic manifold, and let $\me$ be a Riemannian metric on $X$, compatible with the polysymplectic structure. The following are then equivalent:\\
1) $(X,\omega_1,...,\omega_s,\me)$ is an {\em $s$-\ka manifold}.\\
2) $\nabla_X\omega_j = 0$ for all vector fields $X$ and $j = 1,...,s$
\end{teo}
\textit{Proof}
The case $s=1$ is classical, and we therefore omit the proof.\\
$s \geq 2$;\\
Let now $M$ be a smooth manifold of dimension $n(s+1)$, with $s > 1$,  and  let $\omega_{1},...,\omega_{s}$ and $\mathbf{g}$ be as defined in  condition $2)$. Let $p$ be a point of $M$. Pick any standard  polysymplectic coordinate system $x_{i},y^{j}_{i}$,$i~=~1,...,n$,  $j~=~1,...,s$ centered at $p$, defined on a neighborhood $\mathcal{U}$ of $p$ and such that:\\
1) $\forall j ~ ~\omega_{j}~=~\sum_{i}dx_{i}\wedge dy^{j}_{i}$,\\
2) $\mathbf{g}_{p} = \sum_i dx_idx_i~+~\sum_{ij}dy_i^jdy_i^j$\\
i.e. such that the induced coframe on $T_pM$ is orthonormal. Such a  coordinate system exists from the definition of almost $s$-\ka  manifold and from Theorem ~\ref{teo-multisympnorm}. From the fact that $\nabla\omega_j = 0$ for all $j$, we deduce that  parallel transport preserves the polysymplectic structure, and  therefore it must preserve also the standard subspaces associated to  it, among  which are the
\[<\frac{\partial}{\partial  y^{1}_{1}},\ldots,\frac{\partial}{\partial y^{1}_{n}}>,\ldots,<\frac{\partial}{\partial  y^{s}_{1}},\ldots,\frac{\partial}{\partial y^{s}_{1}}>\] 
From this we deduce that for any vector field $X$
\[\nabla_{X}\frac{\partial}{\partial  y^{1}_{i}} = \sum_{l}dy^{1}_{l}\left(\nabla_{X}\frac{\partial}{\partial  y^{1}_{i}}\right)\frac{\partial}{\partial  y^{1}_{l}},\ldots, \nabla_{X}\frac{\partial}{\partial  y^{s}_{i}} = \sum_{l}dy^{s}_{l}\left(\nabla_{X}\frac{\partial}{\partial  y^{s}_{i}}\right)\frac{\partial}{\partial  y^{s}_{l}}\]
As a consequence,  $\nabla_{\frac{\partial}{\partial x_i}}dx_l = -\sum_m  \Gamma_{im}^ldx_m$, where $\Gamma_{im}^l =  dx_l\left(\nabla_{\frac{\partial}{\partial  x_i}}\frac{\partial}{\partial x_m}\right)$ are the usual Christoffel  symbols. We will use the index notation $1,\ldots,n,$ $(11),\ldots,(ns)$  to indicate the $n(s+1)$ indices for the coordinates $x_{i},y^{j}_{i}$,$i~=~1,...,n$,  $j~=~1,...,s$. The above considerations then amount to the fact that  $\Gamma^{(ij)}_{\alpha m} = 0$ for any index $\alpha$, any numbers  $i,m$ in the set $\{1,\ldots,n\}$ and any number $j$ in the set $\{1,\ldots,s\}$. Consider now a coordinate change of the form
\[\tilde{x}_i = x_i + \sum_{mp}b^i_{mp}x_mx_p,~ ~ ~\tilde{y}^j_i =  y^j_i(x_1,...,x_n,y^j_1,...,y^j_n)\]
where the functions $\tilde{y}^j_i$ are determined according to the  Theorem of Carath\'{e}odory-Jacobi-Lie (~\cite{LM} Theorem 13.4 Page 136), so that 
\[\omega_j =  \sum_i d\tilde{x}_i\wedge d\tilde{y}^j_i, ~ ~\tilde{y}^j_i(0,...,0)  = 0\] 
Note that it is crucial that the functions  $\tilde{x}_{i}$ are in involution with respect to the Poisson  structures associated (in the respective  $x_{1},\ldots,x_{n},y^{j}_{1},\ldots,y^{j}_{n}$ spaces) to the various  symplectic forms $\omega_{1},\ldots,\omega_{s}$. In view of the  previous considerations, we see that also in the new coordinates we  have $\nabla_{\frac{\partial}{\partial  \tilde{x}_m}}\frac{\partial}{\partial  \tilde{x}_p} = \sum_p\tilde{\Gamma}^l_{mp}\frac{\partial}{\partial  \tilde{x}_l}$, if the $\tilde{\Gamma}$ are the Christoffel symbols  in the new coordinates, and moreover\\
$\nabla_{\frac{\partial}{\partial  \tilde{x}_m}}d\tilde{x}_l =  \nabla_{\frac{\partial}{\partial  x_m}}d\tilde{x}_l ~+~O(1) =$\\
$\nabla_{\frac{\partial}{\partial  x_m}}\left(dx_l + \sum_{ip}b^l_{ip}x_idx_p\right) ~+~O(1) =  \sum_p\left(-\Gamma^l_{mp} +  b^l_{mp}\right)d\tilde{x}_p ~+~O(1)$.\\  
As it is also the case that $\nabla_{\frac{\partial}{\partial  \tilde{x}_m}}d\tilde{x}_l = -\sum_p\tilde{\Gamma}^l_{mp}d\tilde{x}_p$, if we choose $b^l_{mp} = \Gamma^l_{mp}(0)$ (which we can do as the connection is torsion-free),  we see that the symbols $\tilde{\Gamma}^l_{mp}$ in the new coordinate  system vanish at the origin. For simplicity, we will indicate  the new coordinates with $x_{i},y^{j}_{i}$, and the Christoffel  symbols associated to them with $\Gamma$, dropping the tilde  everywhere.  We know also that for any index $\alpha$, and indicating with $(~ ~)_0$ the evaluation of a form at $0$,\\
$0 = \left(\nabla_{\alpha}\omega_{j}\right)_{0} = \left(\nabla_{\alpha}\sum_{i}dx_{i}\wedge  dy^{j}_{i}\right)_{0} =   \sum_{i}\left(dx_{i}\wedge (\nabla_{\alpha}dy^{j}_{i})\right)_{0} =$.\\
$-\sum_{i}\left(dx_{i}\wedge (\sum_{m k}\Gamma^{(ij)}_{\alpha (m k)}(0)dy^{k}_{m}~+~ \sum_{m}\Gamma^{(ij)}_{\alpha m}(0)dx_{m})\right)_{0}$\\
From this we deduce that $\Gamma^{(ij)}_{\alpha (m k)}(0) = 0$  and $\Gamma^{(ij)}_{\alpha m}(0) = \Gamma^{(mj)}_{\alpha  i}(0)$ for all $i,j,k,m,\alpha$. We consider therefore the change of  coordinates
\[\tilde{y}^{j}_i = y^{j}_i + \sum_{mp}\Gamma^{(ij)}_{mp}(0)x_mx_p,~ ~  ~\tilde{x}_i = x_{i}\]
In the new coordinates we have
\[\sum_{i}dx_{i}\wedge d\tilde{y}^{j}_{i} = \sum_{i}dx_{i} \wedge (dy^{j}_i + \sum_{mp}\Gamma^{(ij)}_{mp}(0)x_mdx_p) =  \omega_{j},\]
as we showed before that $\Gamma^{(ij)}_{mp}(0) = \Gamma^{(pj)}_{mi}(0)$. All the equations for the  Christoffel symbols that we have deduced so far still hold, because  we did not make any assumption on the $y^{j}_{i}$ when we obtained  them, apart from the fact that we were in polysymplectic coordinates. Moreover, we have that 
\[\left(\nabla_{\frac{\partial}{\partial x_l}}d\tilde{y}^{j}_{i}\right)_{0} =  \left(\nabla_{\frac{\partial}{\partial x_l}}dy^{j}_i +  \sum_{mp}\Gamma^{(ij)}_{mp}(0)\nabla_{\frac{\partial}{\partial x_l}}(x_mdx_p)\right)_{0} = \]
\[\left(-\sum_{m}\Gamma^{(ij)}_{lm}(0)dx_{m} +  \sum_{p}\Gamma^{(ij)}_{lp}(0)dx_{p}\right)_{0} = 0\]
From the previous equation, the symmetry of the Christoffel symbols  coming from the fact that the connection is torsion-free, and the  vanishing properties proved above, we see that all the  Christoffel symbols  vanish at $0$.\\
We know from the compatibility of the polysymplectic structure with  the metric that there is a linear change of coordinates which sends  the given coframe at $0$ to an orthonormal (but still polysymplectic)  one. It follows that the same linear change, applied to the  functions $x_{i},y^{j}_{i}$ will preserve the polysymplectic  property, and will make the  coframe at $0$ orthonormal. Moreover, will not  disrupt  the vanishing property (at $0$) of the Christoffel symbols.\\
On the other hand, from the vanishing at the origin of all the Christoffel symbols (and  the fact that the coordinate coframe at $0$ is orthonormal) it is  straightforward to deduce that  $\mathbf{g}~=~\sum_{i}dx_{i}\otimes dx_{i}~+~\sum_{i,j}dy^{j}_{i} \otimes dy^{j}_{i}~+~\mathbf{O}(2)$.
\qed
\begin{pro}
If $\mathbf{X}$ is $2$-\ka and $t\in\R^+$, then also $\alpha_t(\mathbf{X})$, $\beta_t(\mathbf{X})$ and $\lambda_t(\mathbf{X})$ are $2$-\ka.
\end{pro}
\textit{Proof}\\
The statement can be proved locally, where it is clear, using any one of the characterizations of $2$-\ka manifolds.
\qed

From Remark ~\ref{rotatesd} and Proposition ~\ref{covtosd} we obtain the following two remarks:
\begin{rmk}
Let $\mathbf{X} = (X,\omega_1,\omega_2,\me)$ be an almost $2$-\ka manifold with a dualizing form $\omega_D$. If $\mathbf{X}$ is $2$-\ka (i.e. if $\omega_1,\omega_2$ are covariant constant), then $\omega_D$ is covariant constant. In particular, $(X,\omega_1,\omega_2,\me)$ is weakly self dual.
\end{rmk}
\begin{rmk}
On a $2$-\ka manifold there are three different $2$-\ka structures, which are all, from the previous remark, weakly self-dual.
\end{rmk}
We expect that $2$-\ka manifolds will show up as limits of self-dual ones, at limit point of the moduli space where there is some control on the diameter of the manifold. For this reason we expect that the representation on cohomology of $2$-\ka manifolds described in the next section should  be preserved on the monodromy invariant part of the cohomology near well-behaved singularities of almost $2$-\ka manifolds.  
\begin{rmk}
The following question arises naturally from the above results: given a self-dual manifold $(X,\omega_1,\omega_1,\me)$, are there obstructions to deforming $\me$ to a new metric $\mathbf{h}$ for which $(X,\omega_1,\omega_2,\mathbf{h})$ is $2$-\ka?
\end{rmk}
\section{A representation of $\mathbf{sl}(4,\rnum)$}
\label{sec:reps}
In this section we define  a family of operators  (together with their adjoints and associated commutators) which  generalize to $s\geq 1$ the standard Lefschetz operator of \ka  manifolds. Throughout the first part of this section, we assume fixed an almost 2-\ka manifold $(X,\omega_1,\omega_2,\me)$.
\begin{dfn}
\label{llambdah}
The operators $\lop_0,\lop_1,\lop_2$ acting on $\Lambda^*T^*X$, for an almost 2-\ka manifold $(X,\omega_{1},\omega_{2},\mathbf{g})$  are defined  as 
\[\lop_0(\alpha) = \omega_D\wedge\alpha,~ \lop_{1}(\alpha)~=~\omega_{1}\wedge\alpha,~ \lop_{2}(\alpha)~=~\omega_{2}\wedge\alpha\]
\end{dfn}
In the first part of this section we will prove the following:
\begin{teo}
\label{slonforms}
The operators $\lop_0,\lop_1,\lop_2,\lop_0^*,\lop_1^*,\lop_2^*$ generate a Lie algebra naturally isomorphic to $\mathbf{sl}(4,\rnum)$ acting on the bundle $\Lambda^*T^*X$
\end{teo}
\begin{rmk}
\label{rmk-adjoint}
To define the adjoint $\lop_i^*$ to $\lop_i$, we simply used the {\em pointwise}  definition $\forall \alpha,\beta~\in~\Lambda^*T^*_p(X)$~ ~ $\me(\lop_{i}\alpha,\beta)~=~\me(\alpha,\lop_i^*\beta)$. 
\end{rmk}
Before going into the proof, let us remark that the same methods that we will use would show that on an almost $s$-\ka manifold the similarly defined operators generate the real Lie algebra associated to $D_{s+1}$ on the fibres (and the smooth global sections) of $\Lambda^*T^*X$. Note also that the methods of proof of this section are similar to the ones that are used to show that on the complex cotangent bundle of a \ka manifold you have a $\mathbf{sl}(2)$-action. 
\begin{dfn}
1) For $k\in\{0,1,2\}$ the operators $\eik,\iik$ on  $\Lambda^*T^*_pX$, for a  orthonormal standard polysymplectic coordinate coframe $dx_1,...,dx_n,dy^1_1,...,dy^2_n$  at $p$, are  defined as
\[\ei(\alpha)~=~dx_{i}\wedge\alpha,~  \eik(\alpha)~=~dy^{k}_{i}\wedge\alpha~for~k\in\{1,2\},\]
\[\ii~=~\frac{\partial}{\partial x_i}\contr\alpha,~ \iik(\alpha)~=~\frac{\partial}{\partial  y^{k}_{i}}\contr\alpha~for~k\in\{1,2\}\]
2) The operators $\lop_{\alpha\beta}$ on $\Lambda^*T^*_pX$ with $\alpha,\beta\in\{0,1,2,\bar{0},\bar{1},\bar{2}\}$ are defined as
\[\lop_{\alpha\beta}(\phi) = \sum_i\eial\eibe(\phi)~ for~ \alpha,\beta\in\{0,1,2,\bar{0},\bar{1},\bar{2}\}\]
\end{dfn}
\begin{lem}
The operators $\lop_{\alpha\beta}$ are well defined independently of a choice of a basis, and are therefore also well defined as operators acting on the bundle $\Lambda^*T^*X$
\end{lem}
\textit{Proof}
The almost $2$-\ka structure determines canonical isomorphisms of subspaces
\[<dx_1,...,dx_n> \cong <dy^1_1,...,dy^1_n> \cong <dy^2_1,...,dy^2_n>\]
and their duals. The operators $\lop_{\alpha\beta}$ can be interpreted as canonical symplectic forms on spaces of the form $V\times V^\vee$ using these identifications, and are therefore well defined independently of a choice of an orthonormal standard polysymplectic basis.
\qed
\begin{rmk}
1) For any differential form $\phi\in\Omega^*{X}$, we have
\[\lop_{01}(\phi) = \lop_1(\phi) = \omega_1\wedge\phi,~\lop_{02}(\phi) = \lop_2(\phi) = \omega_2\wedge\phi,~ \lop_{12}(\phi) = \lop_0(\phi) = \omega_D\wedge\phi\]
2) \[\lop_{\alpha\beta}^* = \lop_{\bar{\beta}\bar{\alpha}}~ ~for~\alpha,\beta\in\{0,1,2,\bar{0},\bar{1},\bar{2}\}\]
\end{rmk}
The reasoning in the proofs that follow in this section is very similar to the one  that applies to \ka manifolds, used for example in  ~\cite[Pages 106-114]{GH}.
\begin{lem}
\label{rel-among-e-i}
The following relations hold among the operators  $\eial$:\\
1)  $ \eial\ejbe~=~-\ejbe\eial ,~ ~ for~(i,\alpha)~\not=~(j,\bar{\beta})$ \\
2)  $ \eial\eibaral~+~\eibaral\eial~=~Id$ for all $i\in\{1,...,n\},\alpha\in\{0,1,2,\bar{0},\bar{1},\bar{2}\}$
\end{lem}
{\it Proof}~
These identities are easily verified, using the anti-commutativity  property of the wedge product.
\qed
\begin{lem}
\label{commutations}
For $\{\alpha,\beta,\gamma\}\subset\{0,1,2,\bar{0},\bar{1},\bar{2}\}$,\\
1) $\lop_{\alpha\beta} = -\lop_{\beta\alpha}$ when $\alpha \not= \bar{\beta}$\\
2)  $[\lop_{\alpha\beta},\lop_{\bar{\beta}\bar{\alpha}}]~~=~\lop_{\beta\bar{\beta}} -~\lop_{\bar{\alpha}\alpha}$ when $\alpha\not= \beta$\\
3) $[\lop_{\alpha\gamma},\lop_{\bar{\gamma}\beta}]~=~\lop_{\alpha\beta}$ when  $\alpha \not= \gamma,\beta \not= \bar{\gamma}$\\
4) $[\lop_{\alpha\beta},\lop_{\gamma,\delta}] = 0$ when $\{\alpha,\beta\}\cap\{\bar{\gamma},\bar{\delta}\} = \emptyset$\\
5) $[\lop_{\alpha\gamma},\lop_{\gamma\beta}] = 0$ when $\alpha \not= \bar{\gamma}, \beta \not= \bar{\gamma}$
\end{lem}
\textit{Proof}~\\
We can restrict to $\Lambda^*T^*_pX$ for $p\in X$, and use the operators $\eial$ with respect to some (any) orthonormal standard polysymplectic basis.\\
1) 
Immediate from Lemma ~\ref{rel-among-e-i}.\\
2) 
$[\lop_{\alpha\beta},\lop_{\bar{\beta}\bar{\alpha}}]~=~\left(\sum_{i}\eial\eibe\right)\left(\sum_{j}\ejbarbe\ejbaral\right)~-~ \left(\sum_{j}\ejbarbe\ejbaral\right)\left(\sum_{i}\eial\eibe\right)~=$\\
$~=~\sum_{i\not=j}\left(\eial\eibe  \ejbarbe\ejbaral~-~\ejbarbe\ejbaral\eial\eibe\right)~+~ \sum_{i}\left(\eial\eibe\eibarbe\eibaral~-~\eibarbe\eibaral\eial\eibe\right)~~=$\\
$=~\sum_{i}\left(\eial\eibe\eibarbe\eibaral~-~\eibarbe\eibaral\eial\eibe\right)~~~=~ \sum_{i}\left(\eial\eibaral\eibe\eibarbe~-~\eibarbe\eibe\eibaral\eial\right)~= $\\
$=~\sum_{i}\left(\eibe\eibarbe~-~(\eibe\eibarbe~+~\eibarbe\eibe)\eibaral\eial\right)~=~\sum_i\left(\eibe\eibarbe~ -~\eibaral\eial \right)$\\
3) 
$[\lop_{\alpha\gamma},\lop_{\bar{\gamma}\beta}] =  \left(\sum_{i}\eial\eiga\right)\left(\sum_{j}\ejbarga\ejbe\right) - \left(\sum_{j}\ejbarga\ejbe\right)\left(\sum_{i}\eial\eiga\right) =$\\
$= \sum_{i\not=j}\left(\eial\eiga  \ejbarga\ejbe - \ejbarga\ejbe\eial\eiga\right) +  \sum_{i}\left(\eial\eiga  \eibarga\eibe - \eibarga\eibe\eial\eiga\right) =$\\ 
$= \sum_{i}\left(\eial\eiga  \eibarga\eibe - \eibarga\eibe\eial\eiga\right) =
\sum_{i}\eial\eibe\left(\eiga\eibarga+\eibarga\eiga\right) =  \lop_{\alpha\beta}$\\
4)
Immediate from Lemma ~\ref{rel-among-e-i}.\\
5)  
$[\lop_{\alpha\gamma},\lop_{\gamma\beta}] =  \left(\sum_{i}\eial\eiga\right)\left(\sum_{j}\ejga\ejbe\right) - \left(\sum_{j}\ejga\ejbe\right)\left(\sum_{i}\eial\eiga\right) =$\\
$= ~\sum_{i\not=j}\left(\eial\eiga\ejga\ejbe - \ejga\ejbe\eial\eiga\right) +  \sum_{i}\left(\eial\eiga  \eiga\eibe - \eiga\eibe\eial\eiga\right) = 0$\\
\qed

\textit{Proof of Theorem ~\ref{slonforms}} \\
To identify the Lie algebra generated by the $\lop_{\alpha\beta}$ with $\mathbf{sl}(4,\rnum)$ we first writine a Chevalley basis $e_0,e_1,e_2,f_0,f_1,f_2,h_0,h_1,h_2$ for $\mathbf{sl}(4,\rnum)$ satisfying 
\[[e_i,f_j] = \delta_{ij}h_i,~[e_i,h_j] = a_{ij}e_i,~ [f_i,h_j = -a_{ij}f_i\]
(where $(a_{ij})$ is the Cartan matrix for $A_3$) in terms of the $\lop_{\alpha,\beta}$. Let 
\[\begin{array}{lll}
e_0 = \lop_{0\bar{1}} & e_1 = \lop_{1\bar{2}} & e_2 = \lop_{\bar{1}\bar{0}}\\
f_0 = \lop_{\bar{0}1} & f_1 = \lop_{\bar{1}2} & f_2 = \lop_{10} \\
h_0 = \lop_{\bar{0}0} - \lop_{\bar{1}1} & h_1 = \lop_{\bar{1}1} - \lop_{\bar{2}2} & 
h_2 = \lop_{1\bar{1}} - \lop_{\bar{0}0}
\end{array}\]
We then have from Lemma ~\ref{commutations} 2) that $[e_i,f_i] = h_i$ and from 4) $[e_i,f_j] = 0$ for $i\not= j$. Moreover,\\
$[e_0,h_0] = [\lop_{0\bar{1}}, \lop_{\bar{0}0} - \lop_{\bar{1}1}] = 
-[\lop_{\bar{1}0},\lop_{\bar{0}0}] + [\lop_{0\bar{1}},\lop_{1\bar{1}}] = -\lop_{\bar{1}0} + \lop_{0\bar{1}} = 2e_0$\\
$[e_1,h_0] = [\lop_{1\bar{2}}, \lop_{\bar{0}0} - \lop_{\bar{1}1}] = 
[\lop_{\bar{2}1},\lop_{\bar{1}1}] = \lop_{\bar{2}1} = -e_1$\\
$[e_2,h_0] = [\lop_{\bar{1}\bar{0}}, \lop_{\bar{0}0} - \lop_{\bar{1}1}] = 
[\lop_{\bar{0}\bar{1}},\lop_{1\bar{1}}] + [\lop_{\bar{1}\bar{0}},\lop_{0\bar{0}}] = \lop_{\bar{0}\bar{1}} + \lop_{\bar{1}\bar{0}} =  0$\\
$[e_0,h_1] = [\lop_{0\bar{1}}, \lop_{\bar{1}1} - \lop_{\bar{2}2}] = 
-[\lop_{0\bar{1}},\lop_{1\bar{1}}] = -\lop_{0\bar{1}} = -e_0$\\
$[e_1,h_1] = [\lop_{1\bar{2}}, \lop_{\bar{1}1} - \lop_{\bar{2}2}] = 
-[\lop_{\bar{2}1}\lop_{\bar{1}1}] + [\lop_{1\bar{2}},\lop_{2\bar{2}}] = -\lop_{\bar{2}1} + \lop_{1\bar{2}} = 2e_1$\\
$[e_2,h_1] = [\lop_{\bar{1}\bar{0}}, \lop_{\bar{1}1} - \lop_{\bar{2}2}] = 
[\lop_{\bar{0}\bar{1}},\lop_{1\bar{1}}] = \lop_{\bar{0}\bar{1}} = -e_2$\\
$[e_0,h_2] = [\lop_{0\bar{1}}, \lop_{1\bar{1}} - \lop_{\bar{0}0}] = 
 [\lop_{0\bar{1}}, \lop_{1\bar{1}}] + [\lop_{\bar{1}0},\lop_{\bar{0}0}] = \lop_{0\bar{1}} + \lop_{\bar{1}0} = 0$\\
$[e_1,h_2] = [\lop_{1\bar{2}}, \lop_{1\bar{1}} - \lop_{\bar{0}0}] = 
 [\lop_{\bar{2}1},\lop_{\bar{1}1}] ~=~\lop_{\bar{2}1} = -e_1$\\
$[e_2,h_2] = [\lop_{\bar{1}\bar{0}}, \lop_{1\bar{1}} - \lop_{\bar{0}0}] = 
-[\lop_{\bar{0}\bar{1}},\lop_{1\bar{1}}] + [\lop_{\bar{1}\bar{0}},\lop_{0\bar{0}}] = -\lop_{\bar{0}\bar{1}} + \lop_{\bar{1}\bar{0}} = 2e_2$\\
\qed

We now complete the set of  identities which we begun to  describe in Lemma ~\ref{commutations}. These last identities will  allow us to show that we have a representation of the Lie algebra  $\mathbf{sl}(4,\rnum)$ on the cohomology of an s-\ka manifold,  induced by the representation on the space of forms described in  Theorem ~\ref{slonforms}. This will be done showing that the  Laplacian $\Delta_{d}$ commutes with the action of  $\mathbf{sl}(4,\rnum)$. 
\begin{teo}[2-\ka identities]
\label{skaid} ~~ \\
 Let $(X,\omega_{1},\omega_{2},\mathbf{g})$ be an oriented 
s-\ka manifold. Then we have that:\\
1) 
\[[\lop_{hk},d]~=~0~ ~\forall \{h,k\}\subset\{0,1,2\}\]
2)If we define 
$d^{c}_{hk}~:=~[\lop_{hk},d^{*}]$,
we have that\\ 
\[dd^{c}_{hk}~+~d^{c}_{hk}d~=~0~ ~\forall \{h,k\}\subset\{0,1,2\}\]
3) 
\[[\lop_{hk},\Delta_{d}]~=~[\lop_{\bar{k}\bar{h}},\Delta_{d}]~=~0~ ~\forall  \{h,k\}\subset\{0,1,2\}~ ~\]
where $\Delta_{d}$ is the $d$-Laplacian relative to the metric  $\mathbf{g}$ and to the orientation.\\
4) 
\[[\lop_{\alpha\beta},\Delta_d] = 0~ ~\forall\{\alpha,\beta\}\subset\{0,1,2,\bar{0},\bar{1},\bar{2}\}\]
\end{teo}
{\it Proof}\\
1) This equation follows immediately from the fact that $d\omega_1 = d\omega_2 = d\omega_D =  0$.\\
2) If we write down the expression for $d^{c}_{hk}$ in standard 2-\ka coordinates centered  at a point $p\in X$, we see that no derivative of the metric appears. Therefore, when we write down the expression   for $dd^{c}_{hk}~+~d^{c}_{hk}d$, only the first derivatives of the  metric are involved. We skip the details, as they are completely analogous to those of, for example, ~\cite[Pages 111-115]{GH}.\\ It follows, as in the classical case of \ka manifolds, that  to prove the equation it is enough to  reduce to the case of a constant metric.  When the metric is flat, however, the  equation is easily seen to be equivalent (using $1)$) to  $[\lop_{hk},\Delta_d]=0$, which with a flat metric follows immediately from  the fact that $\lop_{hk}$ is constant in flat (orthonormal)  coordinates.\\
3) The second equation is the adjoint of the first. The first one,  once written down explicitely in terms of $d$ and $d^*$, follows  immediately from points $1)-2)$.\\
4) This follows from the previous point, the Jacobi identity and the fact that the Lie algebra of the $\lop_{\alpha,\beta}$ is generated by $\{\lop_{hk}\}\cup\{\lop_{\bar{h}\bar{k}}\}$
\qed
\begin{cor}
Let $(X,\omega_{1},\omega_{2},\me)$ be an oriented  2-\ka  manifold. Then there is a canonical representation of the simple Lie  algebra $\mathbf{sl}(4,\rnum)$ on the space $\mathcal{H}^*(X,\rnum)$ of harmonic forms on $X$
\end{cor}
There is a clear similarity between the representation of $\mathbf{sl}(4,\rnum)$ described in this section and the representations described in ~\cite{LL}. Namely, in both cases one obtains a semi-simple Lie algebra starting from an abelian set of generators, by adding their "$\mathbf{sl}(2)$ adjoints", which still commute among each other. And the space on which these operators act is itself a graded algebra. However, in our case it seems that the representation that you obtain is not a Jordan-Lefschetz module (see ~\cite{LL} for the definition), because, even if separatedly the operators satisfy a form of Lefschetz duality, there does not seem to exist a unique grading associated to the dualities of all of them. In any case, it would be interesting to investigate the connections with the cited  work.\\
\bigskip
~  ~\\
Michele Grassi (grassi@dm.unipi.it)\\
Dipartimento di Matematica, Universit\`{a} di Pisa\\
Via Buonarroti, 2\\ 
56100 Pisa - Italy

\end{document}